\documentclass[letterpaper, 10 pt, conference]{ieeeconf}  
\IEEEoverridecommandlockouts

\usepackage{graphicx}
\usepackage[utf8]{inputenc}
\usepackage[T1]{fontenc}
\usepackage{amsmath}
\usepackage{amsfonts}
\usepackage{amssymb}
\usepackage{mathrsfs}
\usepackage{mathtools}
\usepackage{xcolor}

\newcommand{\norm}[1]{\left\lVert#1\right\rVert}
\newtheorem{assumption}{Assumption}

\newtheorem{remark}{Remark}[section]
\newtheorem{proposition}{Proposition}[section]

\newtheorem{theorem}{Theorem}

\usepackage{cite}

\DeclareMathOperator{\interior}{int}

\title{\LARGE 
    On Exact Solutions to the Linear Bellman Equation
}

\author{David Ohlin$^1$, Richard Pates$^1$ and Murat Arcak$^2$
\thanks{\hspace*{-3mm}This work is partially funded by the Wallenberg AI, Autonomous Systems and Software Program (WASP) funded by the Knut and Alice Wallenberg Foun
dation.}%
\thanks{$^1$ Department of Automatic Control, Lund University, Box 118, SE-221 00 LUND, Sweden. The authors are with the ELLIIT Strategic Research Area at Lund University. (e-mail: {\tt\{david.ohlin, richard.pates\}@control.lth.se}).}%
\thanks{$^2$ Department of Electrical Engineering and Computer Sciences, UC Berkeley, 569 Cory Hall, Berkeley, CA 94720, USA (e-mail: {\tt arcak@berkeley.edu}).}%
}

\begin{document}

\maketitle
\thispagestyle{empty}
\pagestyle{empty}

\begin{abstract}
This paper presents sufficient conditions for optimal control of systems with dynamics given by a linear operator, in order to obtain an explicit solution to the Bellman equation that can be calculated in a distributed fashion. Further, the class of Linearly Solvable MDP is reformulated as a continuous-state optimal control problem. It is shown that this class naturally satisfies the conditions for explicit solution of the Bellman equation, motivating the extension of previous results to semilinear dynamics to account for input nonlinearities. The applicability of the given conditions is illustrated in scenarios with linear and quadratic cost, corresponding to the Stochastic Shortest Path and Linear-Quadratic Regulator problems.
\end{abstract}
\vspace*{-1mm}
\section{Introduction}

\noindent
For what classes of optimal control problems can we expect the existence of efficient algorithms? Leading lights of the field are the $A^*$ \cite{hart68heuristic} {and related heuristic algorithms \cite{barto95rtdp} for deterministic and stochastic versions of the shortest path problem, respectively}. In the Linear-Quadratic Regulator problem (LQR) \cite{kalman60contributions} and Linearly solvable Markov Decision Processes (LDP) \cite{todorov06linearly}, {the special structure of the Bellman equation allows for efficient solution methods}. The objective of this paper is to formalize the connection between these seemingly disparate instances, giving a sufficient condition to guarantee that the Bellman equation can be decoupled and solved explicitly. This decoupling is in turn what enables distributed implementation of algorithms where the optimum preserves the sparsity of the dynamics, as is the case for Stochastic Shortest Path (SSP) and LDP.

The favorable properties of the Bellman equation in the case of linear cost and positive linear dynamics is highlighted in {\cite{blanchini23exact}}, \cite{rantzer22explicit}{, showing that an explicit solution exists under these conditions}. In \cite{ohlin24heuristic}, these results are extended to a class of systems that are shown to be equivalent to a generalized version of SSP{, for which the existence of explicit solutions is previously known (see e.g. \cite{bertsekas91stochastic})}. A general framework for explicit solutions to the Bellman equations in the case of linear dynamics is proposed in \cite{pates24cones}, showing similar results for both linear and quadratic cost based on invariant cones. This methodology is closely connected to the theory of monotone systems developed in~\cite{angeli03monotone}, used in \cite{shen17cones} to generalize previous dissipativity results for positive systems to any invariant cone. The recent work \cite{li25semilinear} characterizes a broader class of dynamics in the case of optimal control with linear cost, by allowing a nonlinear dependence on the input signal. Previous results within the field of Dynamic Programming (DP) are extended to this new class under the assumption of an explicit solution. In the present work, we present sufficient conditions for such a solution based on the dynamics, cost and input constraints of the studied system. 

In \cite{todorov06linearly}, the class of LDP is identified as a subset of MDP with cost based on the Kullback-Liebler distance between an underlying autonomous transition function and the controlled dynamics. {This motivates the more general problem statement in the present work as compared to previous results for linear systems on cones in \cite{pates24cones}.} The key feature of an explicit linear equation for the solution to the Bellman equation is leveraged to pose optimal control of the system as an eigenvalue problem. In order to illustrate the general applicability of our conditions, we reformulate the class of LDP as a continuous-state problem and show that the resulting modified cost and input constraints together with the dynamics share the properties of SSP and LQR that allow for decoupling of the Bellman equation. 

The main contributions of this work are summarized as follows, in order of presentation:

\begin{itemize}
    \item  {Formulation of a unifying optimal control problem} for linear operators on invariant cones.
    \item Statement of a sufficient condition for efficient solution of the Bellman equation.
    \item Reformulation of Linearly solvable MDP (LDP) as a continuous-state optimal control problem.
    \item Demonstration of the given condition and derivation of decomposition functions for SSP, LQR and LDP.
\end{itemize}

Section~II introduces the class of analyzed problems. The main result in the form of Theorem~\ref{th:main} is presented in Section~III. In Section~IV, the applicability of Theorem~\ref{th:main} is shown for SSP and LQR. Section~V comprises the reformulation of LDP and shows the applicability of Theorem~\ref{th:main} also in this case. Finally, conclusions are presented in Section~VI.

\subsection{Notation}
\noindent
Let $\mathbb{R}^n_+$ denote the nonnegative orthant. The operators $\exp(\cdot)$ and $\log(\cdot)$ are used to represent the element-wise exponential and logarithm of a vector or matrix. Further, $\odot$ and $\oslash$ are used to express element-wise multiplication and division between objects of {the same dimension}. The operator $\textnormal{diag}(\cdot)$ extracts the diagonal when applied to a matrix, producing a vector. When applied to a vector, the output is a matrix with diagonal given by that vector.

Consider a Hilbert space $\mathcal{H}$. {Given a linear operator $\mathcal{A}$ on $\mathcal{H}$, let $\mathcal{A}^*$ denote its adjoint.} A set ${\mathcal{X}}\subseteq\mathcal{H}$ is a cone if, for any element $x\in{\mathcal{X}}$ and any positive scalar $a$, it holds that $ax\in{\mathcal{X}}$. The cone ${\mathcal{X}}$ is said to be \textit{proper} if it is {pointed (${0\in\mathcal{X}}$ and for any ${x\in\mathcal{X}}$ such that ${-x\in\mathcal{X}}$, ${x=0}$), closed, convex and solid}. Let the \textit{dual cone} of ${\mathcal{X}}$ be defined as
\begin{equation*}
    {\mathcal{X}}^* = \{y: \langle y,x \rangle\ge0 \;\; \forall x\in{\mathcal{X}}\}
\end{equation*}
Here, $\langle\;\cdot\;,\;\cdot\;\rangle$ denotes an inner product on the underlying Hilbert space. It follows that the dual of a proper cone is also proper. {This work deals with the class of \textit{cone linear absolute norms}, as defined in \cite{seidman05conenorm}. Given a proper cone $\mathcal{X}$, these weighted norms take the following form for $w\in\text{int}(\mathcal{X}^*)$:
\begin{equation*}
    \norm{x}_{w} = \inf \{\langle w, u \rangle: -u \preceq x \preceq u \}
\end{equation*}
which, for elements in the primal cone $x\in\mathcal{X}$, simplifies to}
\begin{equation}\label{eq:norm}
    \norm{x}_w := \langle w, x \rangle.
\end{equation} 

Each proper cone ${\mathcal{X}}$ induces a \textit{partial ordering} of the elements contained within it, so that provided $x,y\in{\mathcal{X}}$ we write $x\succeq y$ if $x-y\in{\mathcal{X}}$ and $x\succ y$ if $x-y\in\text{int}({\mathcal{X}})$. {If neither $x\succeq y$ nor $x\preceq y$ hold, $x$ and $y$ are \textit{unordered}. If the ordering induced by $\mathcal{X}$ is such that no two elements of the set are unordered, we say that the set is \textit{totally ordered} by~$\mathcal{X}$.} When applied to a partially ordered set $\mathcal{K}$, let the operation $\min \mathcal{K}$ return {the set of minimal elements}, meaning any element that is either less than ($\preceq$) or unordered in relation to all other elements of $\mathcal{K}$. In the cases treated in the given examples, the ordering relation will typically be the element-wise comparison $x \le y$ induced by $\mathbb{R}^n_+$, or the Löwner order $A\preceq B$ induced by the semidefinite cone. The {requirement on the dynamics of the problem \eqref{eq:optprob}} are, however, general enough to allow for other {invariant} cones, including restrictions of the dual cone in the aforementioned examples, resulting in more exotic orderings. 

\section{Problem setup}
\noindent
{Consider a bounded linear operator $\mathcal{A}_P$ on a finite dimensional Hilbert space $\mathcal{H}$. 
Let $\mathcal{A}_P$ be continuously parameterized by ${P\in\mathcal{P}}$, where $\mathcal{P}$ is compact. Let ${\mathcal{X}\subset{}\mathcal{H}}$ be a proper cone, and $\mathcal{A}_P$ and $h(P)$ be such that for each ${P\in\mathcal{P}}$, ${\mathcal{A}_P:\mathcal{X}\rightarrow{}\mathcal{X}}$ and ${h(P):\mathcal{P}\rightarrow{}\interior{\mathcal{X}^*}}$. We analyze the infinite-horizon optimal control problem}
\begin{equation}
    \label{eq:optprob}
    \begin{aligned}
        {\underset{P\in\mathcal{P}}{\textnormal{Minimize}}} &\;\;\; \sum\limits_{t=0}^{\infty} \norm{x}_{h(P)}\\
        \textnormal{subject to} &\;\;\; {x(t+1) = \mathcal{A}_Px(t), \;\; x(0) = x_0\in\mathcal{X}}.
    \end{aligned}
\end{equation} 
The existence of a partition with specific structure in relation to the cost and dynamics of \eqref{eq:optprob}, such that both are additively separable, is key to enable decomposition of the Bellman equation. This condition is formalized in the following assumption:
\begin{assumption}\label{as:1}
    {There exists a partition $\{\mathcal{P}_i\}_{i=1}^n$, where 
    \begin{equation}\label{eq:part}
        P\in\mathcal{P} \iff P_i \in\mathcal{P}_i \;\;\;\; i = 1,...,n  
    \end{equation}}
    {such that}
    \begin{equation}\label{eq:decomp}
        \mathcal{A}_P = \sum_{i=1}^n\mathcal{A}_{P_i}\;\;\text{and}\;\;h(P) = \sum_{i=1}^n h_i(P_i).
    \end{equation}
    \vspace*{-1mm}
    {Further, for any $\lambda\in\mathcal{X}^*$, let each of the sets}
    {\begin{equation}\label{eq:subset}
        \{h_i(P_i)+\mathcal{A}^*_{P_i}\lambda : P_i\in\mathcal{P}_i\}
    \end{equation}}
    {be totally ordered by $\mathcal{X}^*$.}
\end{assumption}

\begin{remark}
    {While the number of partitions $n$ is in general unrelated to the dimension of the state, choosing these to coincide is often a convenient way to ensure that the sets \eqref{eq:subset} satisfy Assumption \ref{as:1}. In the examples of SSP and LDP, this corresponds to a coordinate-wise decomposition of the optimization problem, while in the case of LQR it is shown below that this condition is fulfilled by a rank-1 decomposition of the Bellman equation.}
\end{remark}
    
\begin{remark}
    Notably, while \eqref{eq:optprob} is limited to linear operators $\mathcal{A}_P$, the requirements on $\mathcal{P}$ and $h(P)$ allow for the modeling of systems with input nonlinearities. This is illustrated in the case of LDP in Section V.
\end{remark}

Next, we introduce two familiar examples of optimal control problems, namely LQR and SSP, to concretize the general formulation \eqref{eq:optprob}. Examples of the partition \eqref{eq:decomp} for these cases are shown in Section IV.

\subsection{The Linear-Quadratic Regulator}

\noindent
Consider the linear dynamics ${y(t+1) = (A+BK)y(t)}$ with ${y\in\mathbb{R}^n}$, for some matrices $A$, $B$, i.e., assuming static feedback ${u = Ky}$ with gain ${K\in\mathbb{R}^{n\times m}}$ chosen freely. {Boundedness of $\mathcal{P}$ is ensured by restricting the optimization to stabilizing gains.}\footnote{{Alternatively, assuming coercivity of $\mathcal{A}_Px$ with regard to $P$ will ensure that the minimum of \eqref{eq:subset} is attained for unbounded $\mathcal{P}$.}} We illustrate here only the case of ${m=n}$ to simplify the exposition. 
In order to obtain a problem on the form \eqref{eq:optprob}, let ${x(t) = y(t)y(t)^\top}$. This gives the expression ${{\mathcal{A}_Kx = (A+BK)x(A+BK)^\top}}$ for the dynamics with the domain $\mathcal{X}$ given by the semidefinite cone. 
Restricting the dynamics matrix $A+BK$ to be invertible yields invariance of $\mathcal{X}$ under $\mathcal{A}_K$. The immediate cost for symmetric weight matrices $Q$ and $R$ is
\begin{align*}
    \norm{x}_{Q+K^\top RK} &= \langle Q+K^\top RK,x \rangle\\
    &= \text{tr}((Q+K^\top RK)^\top x)\\
    &= y^\top (Q+K^\top RK)y
\end{align*}
using the Frobenius inner product on $\mathcal{X}$. This is equivalent to the typical quadratic cost for ${Q+K^\top RK \succ 0}$, which is exactly the condition on the weight in \eqref{eq:norm}. Finding the $K$ that solves \eqref{eq:optprob} is then equivalent to solving the LQR problem.

\subsection{Stochastic Shortest Path}

\noindent
As shown in \cite{ohlin24heuristic}, SSP can be modeled by the optimal control problem
\begin{equation}
    \label{eq:SSP}
    \begin{aligned}
        \textnormal{Minimize} &\;\;\; \sum\limits_{t=0}^{\infty}\left[ s^\top x(t) + r^\top u(t) \right] \;\textnormal{over}\; \{u(t)\}^\infty_{t=0}\\
        \textnormal{subject to} &\;\;\; x(t+1) = Ax(t) + Bu(t)\\
        &\;\;\; u(t) \ge 0, \;\;\; x(0) = x_0\in\mathbb{R}^n_+\\
        & \;\;\; \begin{matrix} \mathbf{1}^\top u_1(t) & \le & E_1^\top x(t) \\ \vdots & & \vdots \\ \mathbf{1}^\top u_n(t) & \le & E_n^\top x(t) \end{matrix}
    \end{aligned}
\end{equation}
Here, the input signal~${u\in\mathbb{R}^m}$ is partitioned into~$n$ subvectors~$u_i$, each containing~$m_i$ {(possibly ${m_i = 0}$)} elements, so that ${m = \sum_{i = 1}^{n} m_i}$. This is a special case of \eqref{eq:optprob}. Given static feedback ${u = Kx}$, we let $\mathcal{A}_Kx = (A+BK)x$ and invariance of ${\mathcal{X} = \mathbb{R}^n_+}$ under the dynamics corresponds to the condition ${A+BK\ge 0}$. This holds given an appropriate choice of the matrix $E$, governing the input constraints in~\eqref{eq:SSP}. The immediate cost is expressed as the weighted 1-norm
\begin{align*}
    \norm{x}_{s+K^\top r} &= (s+K^\top r)^\top |x|\\
    &= (s+K^\top r)^\top x\;\;\text{for}\;\;x\in\mathcal{X}.
\end{align*}
The inclusion ${s+K^\top r\in\text{int}(\mathcal{X^*})}$ can be fulfilled by requiring ${s > 0}$, ${r\ge0}$, guaranteeing observability of the state.

\section{Main Result}

\noindent
The following theorem gives a sufficient condition for the decomposition and explicit solution of the Bellman equation for \eqref{eq:optprob}. 

\begin{theorem}\label{th:main}
Let Assumption \ref{as:1} hold. Then, the following statements are equivalent:
\begin{itemize}
\item[($i$)] {There exists a $P_*\in\mathcal{P}$ such that} the optimal control problem \eqref{eq:optprob}
has a finite cost for every $x_0\in\mathcal{X}$, {that is achieved with $P=P_*$}.
\item[($ii$)] There exists $\lambda \in \mathcal{X}^*$ satisfying the equation
        \begin{equation}
        \label{eq:p}
            \lambda = \sum_{i=1}^n\min_{P_i\in\mathcal{P}_i} h_i(P_i) + \mathcal{A}_{P_i}^*\lambda
        \end{equation}
\end{itemize}
Further, the optimal value of \eqref{eq:optprob} is given by $\norm{x_0}_{\lambda}$, with $\lambda$ solving~\eqref{eq:p}. The optimally controlled dynamics are given by 
    \begin{equation}\label{eq:argmin}
        P_i = \textnormal{arg}\!\!\!\!\!\!\min_{P_i\in\mathcal{P}_i\;\;\;\;\;}\!\!\!\! h_i(P_i)+\mathcal{A}_{P_i}^*\lambda.
    \end{equation}
\end{theorem}

\quad

\begin{proof}
{\textit{(i)}$\implies$\textit{(ii)}: First observe that for any $P\in\mathcal{P}$ 
\begin{equation}\label{eq:costanyp}
\sum\limits_{t=0}^{\infty} \norm{x(t)}_{h(P)}=\norm{\sum_{t=0}^\infty{}x(t)}_{h(P)}.
\end{equation}
Furthermore, if the constraints in \eqref{eq:optprob} are satisfied, then
${\sum_{t=0}^\infty{}x(t)=\sum_{t=0}^\infty{}\mathcal{A}_{P}^tx_0}$. Since ${h(P)\in\interior{\mathcal{X}^*}}$ it then follows that the cost in \eqref{eq:costanyp} is bounded if and only if $\rho(\mathcal{A}_{P})<1$. Since this implies that ${\sum_{t=0}^\infty{}\mathcal{A}_P^t=(I-\mathcal{A}_P)^{-1}}$, it then follows that the cost equals
\begin{equation}\label{eq:costwithstable}
\norm{(I-\mathcal{A}_{P})^{-1}x_0}_{h(P)}=\norm{(I-\mathcal{A}^*_{P})^{-1}h(P)}_{x_0}.
\end{equation}
Putting ${\lambda_*=(I-\mathcal{A}^*_{P_*})^{-1}h(P_*)}$, we then see that
\[
\lambda_* = h(P_*) + \mathcal{A}_{P_*}^*\lambda_*\in\interior{\mathcal{X}^*},
\]
and the optimal cost in \eqref{eq:optprob} is $\norm{\lambda_*}_{x_0}$. We will now show that the right hand side of the above equals ${\min_{P\in\mathcal{P}}h(P) + \mathcal{A}_{P}^*\lambda_*}$. We proceed by contradiction, so assume that
\begin{equation}\label{eq:contradiction}
\lambda_* \neq{} \min_{P\in\mathcal{P}}h(P) + \mathcal{A}_{P}^*\lambda_*.
\end{equation}
By Assumption \ref{as:1} there must exist a $\overline{P}\in\mathcal{P}$ such that
\begin{equation}\label{eq:f1}
\lambda_*\succeq_{\mathcal{X}^*}h(\overline{P})+\mathcal{A}_{\overline{P}}^*\lambda_*.
\end{equation}
This implies that $\rho({\mathcal{A}_{\overline{P}}^*})<1$, and therefore that there exists a $\overline{\lambda}\in\mathcal{X}^*$ such that
\begin{equation}\label{eq:f2}
\overline{\lambda}=h(\overline{P})+\mathcal{A}_{\overline{P}}^*\overline{\lambda}.
\end{equation}
By \eqref{eq:costwithstable}, the cost \eqref{eq:costanyp} achieved by $\overline{P}$ equals $\norm{\overline{\lambda}}_{x_0}$. Combining \eqref{eq:f1} and \eqref{eq:f2} shows that
\[
(I-\mathcal{A}^*_{\overline{P}})(\lambda_*-\overline{\lambda})\succeq_{\mathcal{X}^*}0.
\]
This implies that
\[
\lambda_*-\overline{\lambda}\succeq_{\mathcal{X}^*}0\implies\norm{{\lambda}_*}_{x_0}\geq{}\norm{\overline{\lambda}}_{x_0}.
\]
Therefore either ${\norm{{\lambda}_*}_{x_0}>\norm{\overline{\lambda}}_{x_0}}$ for some ${x_0\in\mathcal{X}}$, or ${\lambda_*=\overline{\lambda}}$. The first case contradicts the fact that $\norm{{\lambda}_*}_{x_0}$ is the optimal cost. However the second implies that \eqref{eq:f1} holds with equality, which contradicts \eqref{eq:contradiction}. Therefore ${\lambda_* = \min_{P\in\mathcal{P}}h(P) + \mathcal{A}_{P}^*\lambda_*}$. Furthermore, by Assumption~\ref{as:1} it holds that
\begin{equation*}
    \min_{P\in\mathcal{P}}h(P) + \mathcal{A}_{P}^*\lambda_* = \sum_{i=n}^n\min_{P_i\in\mathcal{P}_i}h_i(P_i) + \mathcal{A}_{P_i}^*\lambda_*
\end{equation*}
as required for \textit{(ii)}.}

{\textit{(ii)}$\implies$\textit{(i)}: Let $\lambda_*$ be the solution to \eqref{eq:p} with an associated minimizer ${P=P_*}$. Since \eqref{eq:p} additionally implies that $\rho(\mathcal{A}_{P_*})<1$, it follows from \eqref{eq:costwithstable} that the optimal cost in \eqref{eq:optprob} is upper bounded by $\norm{\lambda_*}_{x_0}$. We will now show that this quantity also lower bounds the cost in \eqref{eq:optprob}. To this end, observe that it follows from the left-hand-side of \eqref{eq:costwithstable} that the cost in \eqref{eq:costanyp} is given by the solution to
\begin{equation}\label{eq:primalrelax}
 \begin{aligned}
            \textnormal{Minimize} &\;\;\; \norm{x}_{h(P)} \;\textnormal{over}\,x\in\mathcal{X}\\
        \textnormal{subject to} &\;\;\; (I-\mathcal{A}_{P})x=x_0.
\end{aligned}
\end{equation}
This is conic optimization on standard form with Lagrangian dual
\begin{equation}\label{eq:dualrelax}
\begin{aligned}
            \textnormal{Maximize} &\;\;\; \norm{\lambda}_{x_0}\\
        \textnormal{subject to} &\;\;\;h(P)+\mathcal{A}_{P}^*\lambda - \lambda \in\mathcal{X}^*.
        \end{aligned}
\end{equation}
By weak duality, it follows that the cost for \eqref{eq:optprob} is lower bounded by
\begin{align*}
\nonumber{}\inf_{P\in\mathcal{P}}\sup_{\lambda\in\mathcal{X}^*}&\{\norm{\lambda}_{x_0}:h(P)+\mathcal{A}_{P}^*\lambda - \lambda \in\mathcal{X}^*\}\\
&\geq{}\sup_{\lambda\in\mathcal{X}^*}\inf_{P\in\mathcal{P}}\{\norm{\lambda}_{x_0}:h(P)+\mathcal{A}_{P}^*\lambda - \lambda \in\mathcal{X}^*\}.
\end{align*}
Recalling that $\mathcal{P}$ is assumed to be closed and bounded, Assumption \ref{as:1} implies that the set ${\{h(P)+\mathcal{A}_P^*\lambda : P\in\mathcal{P}\}}$ has a minimum element. Thus 
we conclude that}
\begin{multline*}
    {\inf_{P\in\mathcal{P}}\{\norm{\lambda}_{x_0}:h(P)+\mathcal{A}_{P}^*\lambda - \lambda \in\mathcal{X}^*\}}\\
    {=}\;
{\begin{cases}
{\norm{\lambda}_{x_0}}&{\text{if $\min_{P\in\mathcal{P}}\left(h(P)+\mathcal{A}_P^*\lambda\right)-\lambda\in\mathcal{X}^*$}}\\
{-\infty}&{\text{otherwise.}}
\end{cases}}
\end{multline*}
{Therefore the optimal cost in \eqref{eq:optprob} is lower bounded by
\begin{align*}
            \textnormal{Maximize} &\;\;\; \norm{\lambda}_{x_0} \\
        \textnormal{subject to} &\;\;\; \min_{P\in\mathcal{P}}\left(h(P)+\mathcal{A}_{P}^*\lambda\right) - \lambda \in\mathcal{X}^*.
        \end{align*}
The solution to \eqref{eq:p} gives a feasible point for the above. Hence ${\norm{\lambda_*}_{x_0} = \norm{x_0}_{\lambda_*}}$ lower bounds the solution to \eqref{eq:optprob}, which completes the proof.}
\end{proof}

\begin{remark}
    When {neither $h_i$ nor $\mathcal{A}_{P_i}^*\lambda$ are strictly $\mathcal{X}^*$-convex \cite{Boyd_Vandenberghe_2004}, the minimizing argument $P_*$ may not be unique.} The above result then holds for any choice $P$ attaining the minimum, and a convenient minimal argument can be selected. {As a consequence of Assumption~\ref{as:1}, the resulting parameter $\lambda$ is still uniquely determined}. In the case of SSP, the optimal cost may in some cases be attained by any linear combination of two actions. 
\end{remark}
\begin{remark}
    In the formulation of the optimal control problem in \eqref{eq:optprob}, the form of the feedback (e.g. static state feedback) is implicitly defined through the optimization variable $P$. Such a formulation is often suitable in applications where practical considerations dictate the form of the feedback used. However, it can be shown in the examples presented that the obtained controls are in fact the best possible, not just the best of the assumed form. This can be established by noting that in these cases (and all others the authors envision), the condition in Theorem~\ref{th:main}\textit{(ii)} gives a solution to the Bellman equation (to 
    make this connection precise, compare the conditions presented here to those for policy iteration, and also to the approaches in \cite{li25semilinear}).
\end{remark}

\section{Examples}\label{sec:examples}

\noindent
This section seeks to contextualize the condition in {Assumption \ref{as:1}} by showing the form it takes in the examples above. 

\subsection{The Linear-Quadratic Regulator}

\noindent
The quadratic cost of LQR means that the expression for optimal static feedback is in general dense, hampering efficient distribution of the calculation. Nevertheless, we illustrate here that the decomposition \eqref{eq:decomp} exists, resulting in~\eqref{eq:p} taking the form of a Riccati equation reformulated as a minimization over the static feedback gain.

The Bellman equation for infinite-horizon optimal control takes the form
\begin{align*}
    \norm{x}_\lambda &= \min_{K} \norm{x}_{Q+K^\top RK} + \norm{\mathcal{A}_Kx}_{\lambda}\\
    &= \min_{K} \langle Q + K^\top RK + (A+BK)^\top\lambda(A+BK),x\rangle
\end{align*}
where the optimal cost is given by $J^*(x) = \norm{x}_{\lambda^*}$ for some $\lambda^*\succ 0$. This is equivalent to the minimum
\begin{equation}\label{eq:minK}
    \lambda = Q + \min_{K} K^\top RK + (A+BK)^\top \lambda (A+BK)
\end{equation}
where $\min\limits_{K}$ is taken with respect to the partial order defined by the semidefinite cone. Rewriting \eqref{eq:minK}, we get
\begin{equation*}
    \lambda = Q + A^\top\lambda A + \min_{K} \Psi
\end{equation*}
with
\begin{equation*}
    \Psi = K^\top(R+B^\top\lambda B)K + A^\top \lambda BK + K^\top B^\top\lambda A.
\end{equation*}
Focusing on this term we use the Cholesky decomposition of the factor 
\begin{equation*}
    R + B^\top\lambda B = LL^\top
\end{equation*}
to make the change of variables ${\hat{K} = L^\top K}$. Letting ${M = L^{-1}B^\top\lambda A}$ we get the rank-1 decomposition
\begin{align*}
    \Psi &= \hat{K}^\top\hat{K} + M^\top\hat{K} + \hat{K}^\top M\\
    &= \sum_{i=1}^n \hat{k}_{i}\hat{k}_i^\top + m_i \hat{k}_i^\top + \hat{k}_i m_i^\top
\end{align*}
where
\begin{equation*}
    \hat{K} = \begin{bmatrix}
        \hat{k}_1 & \cdots & \hat{k}_n
    \end{bmatrix}^\top,\;\;M = \begin{bmatrix}
        m_1 & \cdots & m_n
    \end{bmatrix}^\top
\end{equation*}
Here, the rows $m_{i}$, by the construction of $L$, depend linearly on $\lambda$. As discussed above, ${\mathcal{A}^*_K\lambda = (A+BK)^\top\lambda (A+BK)}$ is positive definite for any $K$, so no additional constraint on the gain is necessary. The minimization of $\Psi$ can be decomposed as in \eqref{eq:decomp}, giving the following expression for the optimal cost
\begin{equation*}
    \lambda = Q + A^\top\lambda A + \sum_{i=1}^n \min_{\hat{k}_i} \hat{k}_{i}\hat{k}_i^\top + m_i \hat{k}_i^\top + \hat{k}_i m_i^\top,
\end{equation*}
fulfilling the requirements of Theorem \ref{th:main}. Each minimization has solution ${\hat{k}_i = -m_i}$ and summing the resulting rank-1 matrices $-m_im_i^\top$ yields the familiar algebraic Riccati equation.

\subsection{Stochastic Shortest Path}

\noindent
Let~$C$ be the binary matrix 
\begin{equation*}
    C = \begin{bmatrix}
        \mathbf{1}_{m_1}^\top &  &  & \\
         & \mathbf{1}_{m_2}^\top & & \\
        & & \ddots & \\
        & & & \mathbf{1}_{m_n}^\top
    \end{bmatrix}
\end{equation*}
where $\mathbf{1}_{m_i}$ is the dimension-$m_i$ vector of all ones and all elements off the block diagonal are zero. The constraints of \eqref{eq:SSP} can be rewritten as
\begin{align*}
    \begin{bmatrix}
        C \\ -I
    \end{bmatrix}u&\le
    \begin{bmatrix}
        E \\ 0
    \end{bmatrix}x\\
    \iff0 &\le \begin{bmatrix}
        E & -C \\ 0 & -I
    \end{bmatrix}\begin{bmatrix}
        x \\ Kx
    \end{bmatrix}\\
    \iff0 &\le \begin{bmatrix}
        E - CK \\ K
    \end{bmatrix}
\end{align*}
{The constraint set $\mathcal{P}$ takes the form
\begin{align*}
    \mathcal{P} &= \{K:\begin{bmatrix}
        E - CK \\ K
    \end{bmatrix}\ge 0\}.
\end{align*}
The partition \eqref{eq:part} is naturally obtained by the feedback law 
\begin{equation*}
    u = \begin{bmatrix}
        K_1^\top & \cdots K_n^\top
    \end{bmatrix}^\top x
\end{equation*}
where $K_i$ is the local feedback in each state, with constraint}
\begin{equation*}
    \mathcal{P}_i = {\{K_i:\begin{bmatrix}
        E_i - \mathbf{1}^\top K_i \\ K_i
    \end{bmatrix}\ge 0\}}.
\end{equation*}
To find the optimal infinite-horizon cost function we state the Bellman equation with ansatz $J(x) = \norm{x}_{\lambda} = \lambda^\top x$
\begin{align*}
    \lambda^\top x &= \min_{K\in\mathcal{P}} (s+K^\top r)^\top x + \lambda^\top(A+BK)x\\
    \iff \lambda &= s + A^\top\lambda + \sum_{i=1}^n\min_{K_i\in\mathcal{P}_i} K_i^\top(r_i+B_i^\top\lambda)
\end{align*}
where ${B=\begin{bmatrix}B_1 \cdots B_n\end{bmatrix}}$ and ${r^\top\!= \begin{bmatrix} r_1^\top & \cdots & r_n^\top \end{bmatrix}}$ are partitioned like the inputs $u_i$. This gives the equation for the stationary point \eqref{eq:p}.

Next, we analyze a problem formulation unlike the two presented above, in that the dynamics of the original optimal control problem depend nonlinearly on the input. This is due to the interpretation of the control as a modification of an underlying transition function.

\section{Linearly Solvable MDP}

\noindent
The following is a motivating example for the extension of previous results on linear systems, e.g. \cite{pates24cones}, to semilinear dynamics. {Originally presented in \cite{todorov06linearly}, the problem is formulated as a set of discrete states with controlled transition probabilities. Costs are incurred as the Kullback-Liebler divergence of the deviations from an underlying transition probability. The advantages of LDP have been successfully applied to inference methods for stochastic optimal control~\cite{kappen12inference}.} We first state the problem as posed in~\cite{todorov06linearly}, then reformulate it as a continuous-state optimal control problem for any initial distribution ${x_0\ge 0}$. Finally, we show that this system fulfills the requirements of Theorem~\ref{th:main}.

Consider an {optimal control problem with probability distribution ${x\in\mathbb{R}^n_+}$ denoting the likelihood of being located in each state $i$ and} input ${u\in\mathbb{R}^{n^2}}$. Let ${\overline{P} = \left[ \overline{p}_1\cdots\overline{p}_n\right]}$ be a column stochastic matrix representing the autonomous dynamics of the system. The input $u$ is split into $n$ subvectors $u_i\in\mathbb{R}^n$, each governing the transmission from state $i$. {Let $x_i$ denote the $i$th element of the state probability distribution vector $x$.} This controlled transmission is given by 
\begin{equation}\label{eq:input}
    p_i = \overline{p}_i\odot \exp(u_i).
\end{equation}
The resulting dynamics are
\begin{equation}\label{eq:ldpdyn}
    x(t+1) = \sum_{i=1}^n p_i x_i{(t)}
\end{equation}
and the immediate cost
\begin{equation}
    g(x,u) = s^\top x + \sum_{i=1}^n x_ip_i^\top u_i.
\end{equation}
Define the set of allowed inputs
\begin{equation}
    \mathcal{U}(x) = \{u:p_i^\top \mathbf{1} = 1\} \cap \mathcal{E}(x)
\end{equation}
where
\begin{equation*}
    \mathcal{E}(x) = \{u: \overline{p}_{ij} = 0 \implies u_{ij} = 0, x_i = 0 \implies u_i = \mathbf{0}\}
\end{equation*}
is a set of added constraints to uniquely determine the cost when the state is empty or the underlying transition rate is zero. {It is clear that the dependence \eqref{eq:input} is nonlinear in $u$, while the system dynamics \eqref{eq:ldpdyn} are linear in $x$. As such, the dynamics belong to the broader class of semilinear systems, recently studied in \cite{li25semilinear}.} Arranging the input vectors in matrix form $U = \left[ u_1\cdots u_n\right]$ lets us reformulate the above problem, giving dynamics on the form ${x(t+1) = Px(t)}$, with immediate cost
\begin{equation*}
    g(x,U) = s^\top x + \textnormal{diag}(P^\top U)^\top x
\end{equation*}
and input constraints
\begin{equation}\label{eq:u}
    \mathcal{U}(x) = \{U:P^\top\textbf{1} = \mathbf{1}\}\cap\mathcal{E}(x)
\end{equation}
where $P = \overline{P}\odot\exp(U)$. The construction of $P$ as column stochastic implies that the total mass $\mathbf{1}^\top x$ is preserved. This in turn means that the problem \eqref{eq:optprob} only has a solution if at least one element of $s$ is zero, and the corresponding state~$x_i$ is absorbing (${p_{ji} = 0}$ for all $j$). Further, an absorbing zero-cost state must be reachable from all states using the dynamics $\overline{P}$. Stabilizing the system then corresponds to locating all mass in these goal states. By removing the goal states and the corresponding transitions (rows and columns of $P$ and~$U$) we arrive at a reduced system ${x_r(t+1) = P_rx_r(t)}$. In the following, the subscript denoting reduced size is dropped in the interest of clear presentation. 

To accommodate the reduction of the state space, the immediate cost function is modified, becoming
\begin{align}\label{eq:gr}
    \norm{x}_{h(P)} &= s^\top\! x + \textnormal{diag}(P^\top \text{log}(P\oslash \overline{P}))^\top x + \pi^\top x
\end{align}
where
\begin{equation}\label{eq:pi}
    \pi = (I-P^\top)\mathbf{1}\odot\log((I-P^\top)\mathbf{1}\oslash\overline{p}_g)
\end{equation}
and $\overline{p}_g$ is the constant vector of unmodified transition rates from non-goal states to the goal state in the original system. This construction makes the immediate cost equal to that of~\cite{todorov06linearly}, with the final term $\pi^\top x$ representing the cost incurred by transitions to the goal states in the original model. We can simplify the constraint set $\mathcal{E}(x)$, as the cost associated with control in states with ${x_i = 0}$ naturally vanishes in the formulation \eqref{eq:gr}.
Elimination of the goal states thus yields the new constraint set
\begin{equation}
    \mathcal{P} = \{P:P^\top \mathbf{1}\le \mathbf{1},\; P\ge 0\}
\end{equation}
with equality $p_i^\top \mathbf{1} = 1$ for row $i$ if the $i$th element of $\overline{p}_g$ is zero. In words, the reduction makes the rate of transition to goal states in the original model dependent on the transitions between non-goal states as specified by $P$. 

In order to show applicability of Theorem \ref{th:main} we first introduce a proposition concerning the immediate cost:
\begin{proposition}
    The immediate cost
    \begin{equation*}
        \norm{x}_{h(P)} = h(P)^\top x
    \end{equation*}
    with
    \begin{equation}\label{eq:gP}
        h(P) = s + \textnormal{diag}(P^\top \log(P\oslash\overline{P})) + \pi    
    \end{equation}
    is a valid norm for $s>0$.
\end{proposition}
\begin{proof}
    The expression \eqref{eq:gP} is equivalent to \eqref{eq:gr}, reformulated to depend only on the dynamics $P$. The term ${(h(P) - s)^\top x}$ is equal to the KL-divergence between two distributions (see \cite{todorov06linearly}), the controlled and autonomous dynamics, weighted by the state vector, and is thus nonnegative. It follows then from positivity of the dynamics that ${h(P)^\top x\ge 0}$ with equality only in the case ${x = 0}$ as a consequence of~${s>0}$. 
\end{proof}

Having restated the problem, we move on to finding the optimal cost function, illustrating that the prerequisites of Theorem \ref{th:main} are indeed fulfilled and how this enables efficient solution of the Bellman equation. Given the cost function \eqref{eq:gr}, we naturally expect an optimal cost function on the same form. This yields the following infinite-horizon Bellman equation:
\begin{align}
    \lambda^\top x &= \min_{P\in\mathcal{P}} h(P)^\top x + \lambda^\top Px\nonumber\\
    \iff \lambda &= s + \min_{P\in\mathcal{P}} \textnormal{diag}(P^\top\log(P\oslash\overline{P})) + \pi + P^\top\lambda\nonumber\\
     &= s + \sum_{i=1}^n\min_{p_i\in\mathcal{P}_i}p_i^\top(\log(p_i\oslash\overline{p}_i)+\lambda) +\pi_i.\label{eq:logpi}
\end{align}
Here $\pi_i$ denotes the $i$th element of the vector $\pi$, which depends only on the column $p_i$ of the dynamics (see \eqref{eq:pi}). {Equation \eqref{eq:logpi} is precisely on the form \eqref{eq:p}. In the case of dense underlying dynamics ($\overline{P}>0$), the minimizers of \eqref{eq:logpi} are given by}
\begin{equation}
    {p_i = \frac{\overline{p}_i\odot\exp(-\lambda)}{\mathbf{1}^\top \overline{p}_i\odot\exp(-\lambda) + (\overline{p}_g)_i},}
\end{equation}
{where $(\overline{p}_g)_i$ extracts the $i$th element of $\overline{p}_g$.} The partition \eqref{eq:part} of the constraint set is
\begin{equation*}
    \mathcal{P}_i = \{p_i:p_i^\top\mathbf{1}\le1,\; p_i\ge0\}
\end{equation*}
with equality $p_i^\top \mathbf{1} = 1$ if the $i$th element of $\overline{p}_g$ is zero. This fulfills the condition \eqref{eq:decomp}. Minimizing each term of the sum yields the solution to \eqref{eq:optprob} according to Theorem~\eqref{th:main}, which in turn gives the solution to the LDP. 

The optimal value of $\lambda$ in the reformulated problem is given by the solution to
\begin{align}
        \lambda &= s - \log(\overline{P}^\top\!\!\exp(-\lambda)+\overline{p}_g)\\
        \exp(-\lambda) &= \textnormal{diag}(\exp(-s))\overline{P}^\top\!\!\exp(-\lambda)+\overline{p}_g\nonumber
\end{align}
Introducing the changes of variables ${z=\exp(-\lambda)}$ and ${G=\textnormal{diag}(\exp(-s))}$, we arrive at the following affine equation, equivalent to the eigenvalue problem derived in \cite{todorov06linearly}:
\begin{equation}\label{eq:eigprob}
    z = G\overline{P}^\top\!z + \overline{p}_g.
\end{equation}
Here, we note that $\overline{P}$ is substochastic and the eigenvalue corresponding to the Perron vector is strictly less than 1 when a goal node is reachable from all nodes of the original problem. Further, $G\le I$ for $s\ge0$, guaranteeing the existence of a solution to \eqref{eq:eigprob} when \eqref{eq:optprob} has a finite value.

The above reformulation together with the application of Theorem \ref{th:main} shows the close connection between traditional results for optimal control of linear systems \cite{kalman60contributions} and the linear cost achieved for LDP in \cite{todorov06linearly}. In this setting, the affine equation \eqref{eq:eigprob} serves as the analog of the Riccati equation in the context of LQR. Any sparsity in the autonomous dynamics $\overline{P}$ is preserved in the optimal solution, similar to the case of linear cost and dynamics \cite{ohlin24heuristic}.

\section{Conclusions}

\noindent
The presented {sufficient condition} ties together three instances of optimal control problems with the special feature of an explicit equation for the solution of the Bellman equation. In all three cases, this result is shown to hinge on the existence of a decomposition of a function of the problem parameters into a sum of functions which can be independently minimized over the original parameter space. This generalizes previous results on optimal control of systems with linear dynamics to include parameterized linear operators, admitting problems where the dependence on the control signal is nonlinear. The {optimal control formulation} expresses the objective function as a cone linear absolute norm, connecting monotonicity of the dynamics with respect to some cone to the form of the {cost functional}. 

{The general form of the stated optimization problem opens up for extensions beyond the known examples treated herein.} {As demonstrated in \cite{pates24cones}, structured restrictions of the dynamically invariant cones allow for the solution of optimal control problems with additional constraints.} {Furthermore, the interpretation of the linear and quadratic cost cases as control with cost based on the first and second stochastic moment of a linear system indicate the possibility of applying the derived results to higher order stochastic moments, with potential applications in risk-sensitive optimal control and control of multilinear systems.}

\bibliographystyle{IEEEtran}
\bibliography{references}

\end{document}